\documentclass[11pt, amsfonts]{amsart}
\usepackage{amsfonts}
\usepackage{latexsym}
\usepackage{amsmath}
\usepackage{amssymb}
\usepackage{amsthm}
\newtheorem{theorem}{Theorem}[section]
\newtheorem{lemma}[theorem]{Lemma}

\newtheorem{proposition}[theorem]{Proposition}
\newtheorem{corollary}[theorem]{Corollary}
\def\R{{\mathbb R}}
\theoremstyle{definition}

\def\R{{\mathbb R}}
\theoremstyle{remark}

\newtheorem*{note*}{Note}

\makeatletter
\newcommand{\ls}{\leqslant}
\newcommand{\gr}{\geqslant}

\makeatother

\begin{document}
\small

\title{\bf Isomorphic properties of intersection bodies}

\author{A.\ Koldobsky} 
\address{Department of Mathematics\\ University of Missouri\\ Columbia\\ Missouri 65211}
\email{koldobskiya@missouri.edu}

\author{G.\ Paouris} 
\address{Department of Mathematics\\Texas A $\&$ M University\\College Station \\TX 77843- 3368\\}
\email{grigoris\_paouris@yahoo.co.uk}

\author{M.\ Zymonopoulou}\address{ Department of Mathematics\\University of Crete\\ Heraklio\\Crete}
\email{marisa.zym@gmail.com}

\maketitle

\begin{abstract}
\footnotesize We study isomorphic properties of two generalizations of intersection
bodies - the class ${\mathcal{I}}_k^n$ of $k$-intersection bodies in $\mathbb R^n$ and the class ${\mathcal{BP}}_k^n$ of generalized 
$k$-intersection bodies in $\mathbb R^n.$  In particular, we show that all
convex bodies can be in a certain sense approximated by intersection bodies, namely,
if $K$ is any symmetric convex body in $\mathbb R^{n}$ and $1\ls k \ls n-1$ then the outer volume ratio distance from $K$ to the
class ${\mathcal{BP}}_{k}^n$ can be estimated by
$$ {\rm o.v.r.} (K, {\mathcal{BP}}_{k}^{n} ) := \inf\{ \left( \frac{ |C|}{ |K|}\right)^{\frac{1}{n}} : C\in {\mathcal{BP}}_{k}^{n}, K\subseteq C\} 
\ls c \sqrt{ \frac{n}{k} \log{\frac{en}{k}}}  ,$$
where $c>0$ is an absolute constant.

Next we prove that if $K$ is a symmetric convex body in $\mathbb R^n$, $1\ls k\ls n-1$ and its $k$-intersection body
$I_{k}(K)$ exists and is convex, then $$ d_{BM}( I_{k}(K), B_{2}^{n}) \ls c(k), $$ 
where $ c(k)$ is a constant depending only on $k$, $d_{BM}$ is the Banach-Mazur distance, and $B_2^n$ is the unit
Euclidean ball in $\mathbb R^n.$ This generalizes a well-known result of Hensley and Borell. We conclude the paper
with volumetric estimates for $k$-intersection bodies.
\end{abstract}

\bigskip

\section{Introduction} 
Let $K$ be an symmetric star body in $\R^n$. Following Lutwak (\cite{Lut1}), we say that a body
$I(K)$ is the {\it intersection body of $K$}  if the radius of $I(K)$ in every direction
is equal to the volume of the central hyperplane section of $K$ perpendicular 
to this direction, i.e. for every $\xi\in S^{n-1},$
$$\rho_{I(K)} (\xi)= |K\cap \xi^\bot|,$$
where $\xi^\bot$ is the central hyperplane perpendicular to $\xi$ and $|\cdot|$ stands for the volume.
A more general class of {\it intersection bodies} can be defined as the 
closure in the radial metric of the class of intersection bodies of star 
bodies.  

Intersection bodies play an important role in the solution
of the Busemann-Petty problem posed in \cite{BP} in 1956: 
suppose that $K$ and $L$ are origin
symmetric convex bodies in $\R^n$ so that, for every $\xi\in S^{n-1},$
$$|K\cap \xi^\bot|\le |L\cap \xi^\bot|.$$
Does it follow that $|K|\le |L|?$
The problem was completely solved in the end of the 90's, and the answer 
is affirmative if $n\le 4$ and negative if $n\ge 5.$
The solution has appeared as a result of work of many mathematicians
(see  \cite[Ch.8]{Ga3} or \cite[Ch.5]{Kol} for details).
A connection between intersection bodies 
and the Busemann-Petty problem was established by Lutwak (\cite{Lut1}): 
the answer to the Busemann-Petty problem in $\R^n$ is affirmative 
if and only if every symmetric convex body in $\R^n$ is an 
intersection body. 

A more general concept of a $k$-intersection body was introduced 
in \cite{Kl4}, \cite{Kl14}. For an integer $k,\ 1\le k <n$ and star bodies $K, D$ in $\R^n,$
we say that $D$ is the {\it $k$-intersection body of $K$} if for every $(n-k)$-dimensional
subspace $H$ of $\R^n,$ 
$$|D\cap H^\bot|= |K\cap H|.$$
The $k$-intersection body of $K$ is unique, but for $k>1$ it does not always exist.
If the $k$-intersection body of $K$ exists, we denote it by $I_k(K).$
Taking the closure in the radial metric of the class of all $D$'s that appear
as $k$-intersection bodies of star bodies, we define the class of 
{\it $k$-intersection bodies}. 
The class of $k$-intersection bodies is related 
to a certain generalization of the Busemann-Petty problem in the same way as 
intersection bodies are related to the original problem (see \cite{Kl4} for details;
this generalization offers a condition on the volume of sections that allows to
compare the volumes of two bodies in arbitrary dimensions). We denote the class
of $k$-intersection bodies in $\R^n$ by ${\mathcal{I}}_{k}^{n}.$  In \cite{Kl14} the first named 
author also gave a Fourier characterization of ${\mathcal{I}}_{k}^{n}$: $K\in {\mathcal{I}}_{k}^{n}$ if and only 
if $\| \cdot \|_{K}^{-k}$ is a positive definite distribution in $\R^n.$

Another generalization of intersection bodies was introduced by Zhang (\cite{Zh3}).
Let $1\ls k \ls n-1$. The spherical Radon transform $R_{n-k}$ is an operator acting
from $C(S^{n-1})$ to the space $C(G_{n,n-k})$ of continuous functions on the 
Grassmanian:
$$R_{n-k}f(H) = \int_{S^{n-1}\cap H} f(x) dx,\qquad \forall f\in C(S^{n-1}),\ H\in Gr_{n,n-k}.$$
We say that a star body $K$ is a {\it generalized $k$-intersection body} 
if there exists a linear positive functional $\mu$ 
on $R_{n-k}(C(S^{n-1}))$ such that for every $f\in C(S^{n-1})$, 
$$ \int_{S^{n-1}} \|x\|_{K}^{-k} f (x) dx = \mu(R_{n-k}(f)).$$ 
Following the notation of \cite{EMilman1} we denote the class of all generalized $k$-intersection bodies 
in $\R^n$ by ${\mathcal{BP}}_{k}^{n}$. 
A characterization of the class ${\mathcal{BP}}_{k}^{n}$ was obtained by Grinberg and Zhang (\cite{GZ}) 
as a generalization of the corresponding result of Goodey and Weil for the original intersection bodies (\cite{GW}): 
a star body $K$ belongs to the class ${\mathcal{BP}}_{k}^{n}$ if and only if it is the limit 
(in the radial metric) of $k$-radial sums of ellipsoids. The class of generalized intersection bodies 
is related to the so-called lower dimensional Busemann-Petty
problem (LDBP-problem). Suppose that $1\le k < n$ and symmetric convex bodies 
$K,L$ in $\R^n$ satisfy 
$$|K\cap H|\le |L\cap H|,$$
for every $(n-k)$-dimensional subspace $H$ of $\R^n.$ Does it follow
that $|K|\le |L| ?$ It was proved in \cite{Zh3} that the answer to this question is affirmative
if and only if every symmetric convex body in $\R^n$ is a generalized $k$-intersection body.
Using this, Bourgain and Zhang (\cite{BZ}) (see also \cite{RZ}) proved that 
for the dimensions of sections $n-k>3$ the answer to the LDBP-problem is negative. 
Another proof of this result was given later in \cite{Kl14}.
The problem is still open in the cases where the dimension of sections $n-k=2,3.$

In the case $k=1$ the classes ${\mathcal{I}}_{1}^{n}={\mathcal{BP}}_{1}^{n}$ coincide with the 
class of original intersection bodies. Also in the case $k=n-1$ the two classes contain all 
star bodies in $\mathbb R^n.$  It was proved in \cite{Kl11} (see also \cite[Corollary 4.9]{Kol})  
that the class ${\mathcal I}_{n-3}^n$ contains all symmetric convex bodies in $\R^n,$ and this is no 
longer true for the classes  ${\mathcal I}_{k}^n$ with $k<n-3$ (\cite{Kl4} or \cite[Theorem 4.13]{Kol}).
It is known that ${\mathcal{BP}}_{k}^{n} \subseteq {\mathcal{I}}_{k}^{n}$ 
(\cite{Kl14}; see also \cite{EMilman1} or \cite[p.92]{Kol} for simpler proofs). The latter two results
immediately imply the negative answer to the LDBP-problem with the dimension
of sections greater than 3. The answer to the problem with two- or three-dimensional sections would
be positive, if the classes ${\mathcal{BP}}_{k}^{n}$ and ${\mathcal{I}}_{k}^{n}$ with $k=n-2$ or $k=n-3$ were equal. 
However, Milman (\cite{EMilman2}) proved that ${\mathcal{BP}}_{k}^{n}$ is a proper subclass of ${\mathcal{I}}_{k}^{n}$ 
for $2\ls k \ls n-2$. The example of Milman is not convex, so LDBP-problem is still open for two- and 
three-dimensional sections. Another open problem
is whether the classes ${\mathcal{I}}_{k}^{n}$ increase with $k.$ It was proved by Yaskin (\cite{Y}) that for $k_1<k_2$
there exists a symmetric convex body that belongs to ${\mathcal{I}}_{k_2}^{n}$ but not to ${\mathcal{I}}_{k_1}^{n}.$
However, the inclusion ${\mathcal{I}}_{k_1}^{n}\subset {\mathcal{I}}_{k_2}^{n}$ is known only in the case where $k_1$
divides $k_2$ (see \cite{EMilman1}). For more results on these classes of bodies see \cite{KY} 
and references there.

In spite of all these results, the isomorphic properties of intersection bodies are not very well understood.
The first result of this kind was established by Hensley and Borell (see \cite{H}, \cite{Bor}): if $K$ is symmetric and convex,
then the Banach-Mazur distance $ d_{BM}(I(K), B_{2}^{n})\ls c$, where $c>0$ is an absolute 
constant, which means that intersection bodies of convex bodies are isomorphic to ellipsoids (note that if $K$ 
is a symmetric convex body, then the classical result of Busemann (\cite{Bus1}) guaranties that $I(K)$ 
is also symmetric and convex).  Busemann (\cite{Bus3}) also showed an 
``isoperimetric" type inequality: if $K$ is symmetric convex and $|K|= 1$, 
then $ |I(K)| \ls | I(D_{n})|$, where $D_{n}$ is the Euclidean ball with volume 1. 
This result can be extended to a class more general than convex bodies (even for Borel measurable sets) 
(see \cite{Grin}, \cite{Gard}).  However, intersection bodies of convex bodies form only a small part of the class 
of convex intersection bodies. 

As it was proved by the first named author (\cite{Kl15})(using the Fourier characterization of the intersection bodies), 
the unit ball of any finite dimensional subspace of $L_p,\ p\in (0,2]$ is a $k$-intersection body for every $k$, and in particular 
all polar projection bodies (unit balls of subspaces of $L_1$) are intersection bodies. The class of intersection bodies is strictly larger 
than the class of polar of projection bodies; see \cite{Kl9} for examples. A long standing question is if the two classes are ``isomorphic", i.e. 
whether for every intersection body $I$ there exists a projection body $\Pi$, such that $ d_{BM}(I, \Pi^{\circ}) \ls c$, where $c>0$ 
is an absolute constant (this question is related to the 1970 problem of Kwapien (\cite{Kw}) from the Banach space theory
through the connection between intersection bodies and $L_p$-spaces; see \cite{KK2}). 
A closely related result was proved in \cite{KK2}: for any $q\in (0,1)$ and any $k<n$, there exists 
a constant $c(q,k)$ such that for any convex $k$-intersection body $K$
in $\mathbb R^n$ there exists a subspace of $L_q$ whose unit ball $D$ satisfies $d_{BM}(K,D)<c(q,k).$ 
Note that the constant $c(q,k)$ goes to infinity when $q$
tends to 1, and, if not for that, the case $q=1$ would imply the desired result for polar projection bodies.

In this article we prove several isomorphic results for intersection bodies and their generalizations. We have already mentioned the fact
that the class ${\mathcal I}_{n-3}^n$ contains all symmetric convex bodies, but this is no longer the case for the classes ${\mathcal I}_{k}^n$ with $k<n-3.$
We start with a result showing that $k$-intersection bodies with $k<n-3$ are still in some sense dense in the class of all convex bodies.

Let $K$ be a symmetric convex body in $\mathbb R^n$ and let ${\mathcal{A}}$ be a class of star bodies. We define the outer volume ratio distance by
$${\rm o.v.r}(K, {\mathcal{A}}) := \inf \{{\rm o.v.r.}(K, C), \ C\in {\mathcal{A}}\}, $$
where $${\rm o.v.r.}(K,C):= \inf\{ \left(\frac{|T C|}{|K|}\right)^{\frac{1}{n}} : K\subseteq TC, \ T\in GL_{n}\} .$$ 

\begin{theorem}
\noindent Let $K$ be a symmetric convex body in $\mathbb R^n$ and $1\ls k \ls n-1$. Then 
$$ {\rm o.v.r.}( K, {\mathcal{BP}}_{k}^{n}) \ls c\sqrt{\frac{n \log{\frac{en}{k}} }{k}} , \leqno(1.1)$$
where $c>0$ is an absolute constant. 
\end{theorem}
\noindent Recall that ${\mathcal{BP}}_{k}^{n}\subseteq {\mathcal{I}}_{k}^{n},$ so the result also applies to $k$-intersection bodies.

\smallskip

Our second result extends to the class of $k$-intersection bodies the classical result of Hensley and Borell that an intersection body of a convex body is isomorphic 
to an ellipsoid. Here one faces two additional difficulties. First the $k$-intersection body of a convex body does not necessarily exist and, secondly, even 
if it exists it may not be convex. So any result must take into account these two conditions as additional assumptions. We prove

\begin{theorem}
\noindent Let $K$ be a symmetric convex body in $\mathbb R^n$, $1\ls k \ls n-1$ and assume that $I_{k}(K)$ exists and it is convex. Then 
$$  d_{BM} (I_{k}(K), B_{2}^{n}) \ls c(k), \leqno(1.2)$$
where $c(k)$ depends only on $k$. 
\end{theorem}
 
\noindent Finally, we get some estimates for the volume radius of $k$-intersection bodies:

\begin{theorem}
\noindent Let $K$ be a symmetric star-shaped body in $\mathbb R^n$ with $|K|= |B_{2}^{n}|$. Assume that the $k$-intersection body of $K$ exists. Then 
$$ \left( \frac{|I_{k}(K)|}{| I_{k}(B_{2}^{n})|}\right)^{\frac{1}{n}} \gr \frac{L_{B_{2}^{n}}}{L_{K}}, \leqno(1.3) $$
with equality if and only if $K$ is a symmetric ellipsoid. Here $L_{K}$ stands for the isotropic constant of $K$. 

\smallskip

\noindent Moreover, if $I_{k}(K)$ is a convex body, we have that
$$  \left( \frac{|I_{k}(K)|}{| I_{k}(B_{2}^{n})|}\right)^{\frac{1}{n}} \ls c \log{\left( 1+ d_{BM}(I_{k}(K), B_{2}^{n})\right)} \ls c\min\{ \log{n}, k\log{k}\} , \leqno(1.4) $$
where $c>0$ is a universal constant.
\end{theorem}

\smallskip

\noindent The paper is organized as follows: In section 2 we introduce basic definitions and notation. In sections 3, 4 and  5 we give the 
proof of Theorems 1.1, 1.2 and 1.3 respectively. We provide some final remarks in section 6.


\section{Notation and Definitions}

 We work in
${\mathbb R}^n$, which is equipped with a Euclidean structure
$\langle\cdot ,\cdot\rangle $. We denote by $\|\cdot \|_2$ the
corresponding Euclidean norm, and write $B_2^n$ for the Euclidean
unit ball, and $S^{n-1}$ for the unit sphere. Volume is denoted
by $|\cdot |$. We write $\omega_n$ for the volume of $B_2^n$ and
$\sigma $ for the rotationally invariant probability measure on
$S^{n-1}$. We will write $D_{n}$ for the euclidean ball of volume $1$. ($D_{n}:= |B_{2}^{n}|^{-\frac{1}{n}} B_{2}^{n}$). The Grassmannian manifold $G_{n,k}$ of $k$-dimensional
subspaces of ${\mathbb R}^n$ is equipped with the Haar probability
measure $\mu_{n,k}$. We denote $GL_{n}$ the set of linear invertible transformations and 
$SL_{n}$ for the measure preserving linear transformations.

A compact set $K$ in $\R^n$ is called a star body if the origin is an interior point of $K$, every straight
line passing through the origin crosses the boundary of $K$ at exactly two points and the Minkowski
functional of $K$ defined by
$$\|x\|_K= \min\{a\ge 0:\ x\in aK\},\qquad \forall x\in \R^n$$
is a continuous function on $\R^n.$ The radial function of $K$ is given by 
$$\rho_K(x) = \max\{a>0:\ ax\in K.$$
Throughout this article we say that $K$ is symmetric if $K=-K.$
The radial metric on the class of star bodies is defined by
$$d_r(K,L) =\max_{\xi\in S^{n-1}} |\rho_K(\xi)-\rho_L(\xi)|.$$
If $K,L$ are star bodies and $k\in \mathbb{N},$ then the $k$-radial sum of $K$ and $L$
is a star body $D$ whose radial function is equal to
$$\rho_D = \left( \rho_K^k + \rho_L^k\right)^{1/k}.$$
If $K$ is a convex body, we write $h_{K}$ for the support function of $K$:
$$h_K(x) =\max_{\xi\in K} \langle x,\xi\rangle \qquad \forall x\in \R^n.$$
We define the geometric distance between two bodies $K_{1}$ and $K_{2}$, $d_{G}(K_{1},K_{2})$ as the infimum of positive numbers $r$ 
such that there exists some $a>0$ so that
$$ K_{1}\subseteq aK_{2} \subseteq raK_{1} . $$
The Banach-Mazur distance is defined as 
$$ d_{BM}(K_{1}, K_{2}) := \inf_{T\in GL_{n}}d_{G}(K_{1}, TK_{2}) . $$ 

\noindent Let $K$ be a star body in $\mathbb R^n$ of volume $1$. We define the isotropic constant of $K$, $L_{K}$ as
$$ L_{K}^{2} := \frac{ 1}{n} \inf\{ \int_{TK}\|x\|_{2}^{2} dx : T \in SL_{n} \} . $$  

\noindent For two convex bodies K and L in $\R^n$, the covering number of $K$ by $L,$ denoted by $N(K,L),$  is defined as the minimal number 
of translates of $L$ with their centers in K, needed to cover $K.$ 

\noindent The notation $a\simeq b$ means that there exist universal constants $c_{1}, c_{2}>0$ such that $c_{1} a\ls b \ls c_{2} b$.
We refer to the books \cite{Ga3}, \cite{Kol}, \cite{MS}, \cite{Pisier1}, \cite{Sch} for
basic facts from the Brunn-Minkowski theory, the asymptotic
theory of finite dimensional normed spaces and intersection bodies.


\section{On the outer volume ratio of convex bodies with respect to the class ${\mathcal{BP}}_{k}^{n} $}

\noindent The main idea of the proof of Theorem 1.1. is to first cover a given convex body with Euclidean balls and then show that, if the 
number of balls is not too big, one can approximate the union of the balls by a body in ${\mathcal{BP}}_{k}^{n},$ 
where $k$ will be related to the covering number of $K$. 

\noindent We will use the following theorem for covering numbers of Pisier (see \cite{Pisier1}).

\begin{theorem}\label{Pisier}
\noindent Let $K$ a symmetric convex body in $\mathbb R^n$ and $\alpha\in (0,2)$. Then there exist an $T \in SL_{n}$ such that if $K_{1}= TK$ then
\begin{enumerate}
\item $|K_{1}|= |B_{2}^{n}| $,
\item $ \log{N( K_{1}, t B_{2}^{n})} \ls  \frac{c n}{t^{a} (2-\alpha) } $, $ \forall t\gr 1$, where $c>0$ is a universal constant. 
\end{enumerate}
\end{theorem}

\smallskip
\noindent It is known that the position of the body $K,$ constructed in the previous theorem, 
satisfies the reverse Brunn-Minkowski inequality of V.\ Milman (\cite{Mil5}). 
For completeness we provide a proof.

\begin{corollary}
\noindent Let $K$ be a symmetric convex body with the properties of the body $K_1$ from the previous theorem. Then, for $t\gr 1$
$$ |K+ tB_{2}^{n}|^{\frac{1}{n}} \ls c t|K|^{\frac{1}{n}} , \leqno(3.1) $$
where $c>0$ is an absolute constant. 
\end{corollary}

\noindent {\it Proof.}
\noindent Choose $\alpha=1$. Then using Theorem \ref{Pisier}, we have
$$ \frac{  |K+ tB_{2}^{n}|^{\frac{1}{n}}}{  t|K|^{\frac{1}{n}}} = \frac{  |K+ tB_{2}^{n}|^{\frac{1}{n}}}{  t|B_{2}^{n} |^{\frac{1}{n}}}  \ls N\left( K+ tB_{2}^{n} , 2t B_{2}^{n}\right)^{\frac{1}{n}} \ls  $$
$$ N\left( K, tB_{2}^{n}\right)^{\frac{1}{n}} \ls c .$$
$\hfill\Box $ 

\smallskip

\noindent We will need the following elementary

\begin{lemma}
\noindent Let $z\in \mathbb R^{n}$ and $t>0$. Then there exists a centered ellipsoid ${\mathcal{E}}$ such that
$$  z+ t B_{2}^{n} \subseteq {\mathcal{E}} \subseteq co\{ 2 z + 2\sqrt{2}tB_{2}^{n} ,    -2 z + 2\sqrt{2}tB_{2}^{n} \} . \leqno(3.2)$$
\end{lemma}

\noindent {\it Proof.}
\noindent We may assume that $z:= z_{1} e_{1}$. Let $C$ be the ``cylinder" defined as
$$ C:= \{ (s, y)\in \mathbb R^n: |s| \ls z_{1} + t ,\ \|y\|_{2} \ls t \} . $$
Then one can check that 
$$ z+ tB_{2}^{n} \subseteq C \subseteq {\rm co}\{ z+ \sqrt{2}t B_{2}^{n}, -z+ \sqrt{2} t B_{2}^{n} \} =: K \leqno(3.3)$$
Let $Q:=  \{ (s, y)\in \mathbb R^n: |s| \ls 1 ,\ \|y\|_{2} \ls 1 \}  $. Then $C:= T Q$, where $T\in GL_{n}$. Define ${\mathcal{E}}_{1}:= T B_{2}^{n}$. Then 
$$ B_{2}^{n} \subseteq Q \subseteq \sqrt{2} B_{2}^{n} , \ {\rm or} \  {\mathcal{E}}_{1} \subseteq C \subseteq \sqrt{2} {\mathcal{E}}_{1} . \leqno(3.4) $$
Then, by (3.3) we get that
$$ z+ tB_{2}^{n} \subseteq C \subseteq \sqrt{2}  {\mathcal{E}}_{1} \subseteq 2C \subseteq 2 K . $$
We set ${\mathcal{E}}:= \sqrt{2}{\mathcal{E}}_{1}$ and the proof is complete. $\hfill\Box $ 

\smallskip
 
\noindent Note that the class ${\mathcal{BP}}_{k}^{n}$ is closed under $k$-radial sums and that ellipsoids belogn to this class for all $1\ls k \ls n-1$.

\noindent Let $V_{1}, V_{2}$ two star bodies in $\mathbb R^n.$ We define the distance between $V_{1}$ and $V_{2}$ as
$$ d(V_{1}, V_{2}) := \sup_{ x\in \mathbb R^n, x\neq 0 } \big\{ \frac{ \rho_{V_{1}} (x) }{ \rho_{V_{2}}(x) }  ,  \frac{\rho_{V_{2}}(x)}{ \rho_{V_{1}}(x)}  \big\} . \leqno(3.5) $$

\noindent Observe that the definition implies that 
$$ \frac{1}{d} V_{1} \subseteq V_{2} \subseteq d V_{1} , \leqno(3.6) $$
where $d:= d(V_{1}, V_{2})$. 

\begin{lemma}
\noindent  Let $1\ls k \ls n-1$ and $ \log{N} \ls k$. Let $V_{i}, i\ls N,$  be symmetric star bodies in ${\mathcal{BP}}_{k}^{n}$  and let $V= \cup_{i=1}^{N} V_{i}$. 
Then there exists  $C\in {\mathcal{BP}}_{k}^{n}$ such that
$$ d(V, C) \ls e . \leqno(3.7)$$
\end{lemma}

\noindent {\it Proof.}
\noindent Let $\rho_{V_{i}}$ be the corresponding radial functions of $V_{i}$. Then 
$$ \rho_{V} := \max_{i\ls N} \rho_{V_{i}}(x), \ x \in \mathbb R^n\setminus\{0\} . \leqno(3.8) $$
For any $a\in \mathbb R^{N}$ we have that 
$$ \|a\|_{\infty} \ls \|a\|_{p} \ls N^{\frac{1}{k}} \|a\|_{\infty} \ls e \|a\|_{\infty} . \leqno(3.9) $$
Let $C$ be the star body defined by 
$$ \rho_{C}^{k} (x) := \sum_{i=1}^{N} \rho_{V_{i}}^{k} (x) , \ 0\neq x \in \mathbb R^n . \leqno(3.10)  $$
Note that $C\in {\mathcal{BP}}_{k}^{n}$ and by (3.9) we also have that $d(V,C)\ls e$. $\hfill\Box $ 

\smallskip

\noindent We can now give a proof of Theorem 1.1:

\noindent Let  $ \alpha:=2-\frac{1}{\log{e\frac{n}{k}}}$. Note that $ {\rm o.v.r.}(T K, {\mathcal{BP}}_{k}^{n}) = {\rm o.v.r.}(K, {\mathcal{BP}}_{k}^{n} ) $ for any $T \in GL_{n}$. 
So we may assume that $K$ is in the position described in Theorem 3.1. We have that for every $t\gr 1$, there exists $N$, 
such that $ \log{N} \ls \frac{cn}{t^{a}(2-a)} $, and $ z_{1}, \dots , z_{N} \in K$ such that 
$$ K \subseteq \cup_{i=1}^{N} z_{i} + t B_{2}^{n} . $$
Given any $z_{i},$ let $ {\mathcal{E}}_{i}$ as in Lemma 3.3. Then we have that
$$ K \subseteq \cup_{i=1}^{N} {\mathcal{E}}_{i} \subseteq 2 \left( co\{ z_{1}, \dots, z_{n}\} + \sqrt{2} t B_{2}^{n} \right)  \subseteq 2( K + \sqrt{2}t B_{2}^{n}) . $$
 Choose $t$ such that $ k=  \frac{cn}{t^{a}(2-a)} $. Then $\log{N} \ls k$. So, by Lemma 3.4, there exist a $C\in {\mathcal{BP}}_{k}^{n}$ such that $ d(C, V) \ls e$, 
 where $V=  \cup_{i=1}^{N} {\mathcal{E}}_{i}$. So, 
$$ \frac{1}{e} K \subseteq \frac{1}{e} V \subseteq C\subseteq e V \subseteq 2e\left( K + \sqrt{2} t B_{2}^{n}\right) \leqno(3.11) $$
Recall that  $t:= \left( \frac{n}{k(2-a)}\right)^{\frac{1}{a}}\simeq \sqrt{\frac{n}{k}} \sqrt{\log{e\frac{n}{k}}}$. Then (3.11) becomes
$$  \frac{1}{e} K \subseteq C \subseteq 2e\left( K + c  \sqrt{\frac{n}{k}} \sqrt{\log{e\frac{n}{k}}} B_{2}^{n}\right) \leqno(3.12) $$
Moreover, by Corollary 3.2, we have that
$$ |C|^{\frac{1}{n}} \ls 2e  |K + t B_{2}^{n} |^{\frac{1}{n}} \ls c t |K|^{\frac{1}{n}}. $$
Let $C_{1}:= eC$. We have that $C_{1}\in {\mathcal{BP}}_{k}^{n}$, $ K\subseteq C_{1}$ and 
$$ \frac{ |C_{1} |^{\frac{1}{n}} }  {|K|^{\frac{1}{n}} } \ls c^{\prime} t \ls c^{\prime\prime}   \sqrt{\frac{n}{k}} \sqrt{\log{e\frac{n}{k}}} . \leqno(3.13) $$
This finishes the proof.  $\hfill\Box $


\section{Distances for $k$-intersection bodies}

\noindent Let $K$ a star body in $\mathbb R^n$. 
Recall that the $k$-intersection body of $K$ (if it exists) is a star body that satisfies 
$$ |I_{k}(K) \cap F  | = |K\cap F^{\perp}|, \ \forall F\in G_{n,k} . \leqno(4.1) $$

\noindent Note that $I_{k}(K)$ (if it exists) is unique. Also if $t>0$ and $T\in GL_{n},$ one has that 
$$ I_{k}(tK)= t^{\frac{n-k}{k}} I_{k}(K) \ {\rm and}  \  I_{k}(TK)= |{\rm det} T| T^{-\ast}( I_{k}(K))   . \leqno(4.2) $$
The last equality follows from the Fourier characterization of the $k$-intersection bodies (see \cite{Kl14}).

\noindent Let $K$ be  a symmetric convex body of volume $1$ and $p\gr 1$. We define the $L_{p}$-centroid body of $K$ (\cite{LZ},\cite{LYZ},\cite{Pa3}), 
as the symmetric convex body that has support function
$$ h_{Z_{p}(K)}(\theta):= \left(\int_{K} |\langle x, \theta \rangle|^{p} d x\right)^{\frac{1}{p}} . \leqno(4.3)$$

\noindent We will use the following $L_{q}$-version of Rogers-Shephard inequality (\cite{Pa3},\cite{Pa4}):

\begin{proposition}
\noindent Let $K$ be  a symmetric convex body of volume $1$ in $\mathbb R^n$ and $1\ls k \ls n-1$. Then for every $F\in G_{n,k}$,
$$ |K\cap F^{\perp} |^{\frac{1}{k}} |P_{F} Z_{k}(K)|^{\frac{1}{k}} \simeq 1. \leqno(4.4) $$ 
\end{proposition} 

\noindent We will also use the Santal\'o and reverse Santal\'o inequality (\cite{BM2}): If $K$ is a symmetric convex body in $\mathbb R^n$, then 
$$ \left(|K| |K^{\circ}|\right)^{\frac{1}{n}} \simeq \frac{1}{n} . \leqno(4.5) $$

\begin{proposition}
\noindent Let $K$ be a symmetric convex body in $\mathbb R^n$ of volume $1$ and $1\ls k\ls n-1$. Assume that $I_{k}(K)$ exists. Then for $F\in G_{n,k}$, 
$$ |I_{k}(K) \cap F|^{\frac{1}{k}} \simeq | \left( \frac{Z_{k}(K)}{k}\right)^{\circ} \cap F |^{\frac{1}{k}} . \leqno(4.6) $$
\end{proposition}

\noindent {\it Proof.}
\noindent Using the definition (4.1) and equations (4.4), (4.5), we get that 
$$  |I_{k}(K) \cap F|^{\frac{1}{k}} =  |K \cap F^{\perp}|^{\frac{1}{k}} \simeq |P_{F} Z_{k}(K)|^{-\frac{1}{k}} \simeq $$
$$ k|Z_{k}^{\circ}(K) \cap F|^{\frac{1}{k}} =  | \left( \frac{Z_{k}(K)}{k}\right)^{\circ} \cap F |^{\frac{1}{k}} . $$
$\hfill\Box $ 

\noindent The following lemma is a well known application of the Brunn-Minkowski inequality:

\begin{lemma}
\noindent Let $K$ be a symmetric convex body of volume $1$ in $\mathbb R^n$ and let $k\gr 2$. Then 
$$ Z_{2}(K) \subseteq Z_{k}(K)\subseteq ckZ_{2}(K), \leqno(4.7) $$
where $c>0$ is an absolute constant. 
\end{lemma}

\noindent Note that $Z_{2}(K)$ is an ellipsoid. Moreover, since $Z_{p}(TK)= TZ_{p}(K)$ for $T\in SL_{n}$, there exists $T\in SL_{n}$ 
such that $Z_{2}(TK):= L_{K} B_{2}^{n}$. In this case we say that $K$ is isotropic (see \cite{MP}, \cite{Gian} for more information on isotropicity). Note that in the case where $K$ is convex the definition of the isotropic constant $L_{K}$ that we give here is equivalent to the definition given in \S2. (see e.g. \cite{Gian}).

\noindent We have the following

\begin{corollary}
\noindent Let $K$ be an isotropic convex body in $\mathbb R^n$, $1\ls k \ls n-1$ and assume that $I_{k}(K)$ exists. Then for all $F_{1}, F_{2} \in G_{n,k}$, 
$$ \frac{1}{(c_{1}k)^{k}} \ls \frac{ |I_{k}(K)\cap F_{1}|}{|I_{k}(K)\cap F_{2}|} \ls (c_{1}k)^{k}, \leqno(4.8) $$
where $c_{1}> 1$ is an absolute constant. 
\end{corollary}

\noindent {\it Proof.}
\noindent Using equations (4.6) and (4.7) we have that
$$ \frac{  |I_{k}(K)\cap F_{1}|}{|I_{k}(K)\cap F_{2}|} \ls c^{k} \frac{  | \left( \frac{Z_{k}(K)}{k}\right)^{\circ} \cap F_{1} |}{  | \left( \frac{Z_{k}(K)}{k}\right)^{\circ} \cap F_{2} |}  $$
$$ \ls  c^{k} \frac{  | \left( \frac{Z_{2}(K)}{k}\right)^{\circ} \cap F_{1} |}{c_{0}^{k}  | Z_{2}^{\circ}(K)  \cap F_{2} |} \ls \frac{k^{k}c^{k}}{c_{0}^{k}}  \frac{  | Z_{2}^{\circ}(K) \cap F_{1} |}
{  | Z_{2}^{\circ}(K)  \cap F_{2} |} \ls (c_{1}k)^{k} ,$$
since $K$ is isotropic.
We work similary for the left hand side inequality. $\hfill\Box $ 

\smallskip

\noindent We will also need the following

\begin{lemma}
\noindent Let $K$ be a symmetric convex body in $\mathbb R^{m+1}$. Let $r:= r_{K}(e_{m+1})$ and $R:= h_{K}(e_{m+1})$. Then 
$$ \frac{2r}{m+1}|K\cap \mathbb R^{m}|  \ls  |K| \ls 2R |K\cap \mathbb R^{m}| . \leqno(4.9) $$ 
\end{lemma}

\noindent {\it Proof.}
\noindent To estimate the left-hand side observe that $re_{m+1}\in K,$ so
$$ |K| \gr | {\rm co} \{ K\cap \mathbb R^m, r e_{m+1}, - r e_{m+1} \} | = 2\frac{r}{m+1}|K\cap \mathbb R^{m}| . $$
For the right-hand side, observe that the function $f(s):= |K\cap (se_{m+1} + e_{m+1}^{\perp})|$ is even and $log$-concave by the Brunn-Minkowski inequality, therefore  attains its maximum at $0$. 
Then by Fubini's theorem,
$$ |K| = 2 \int_{0}^{R} f(s) d s \ls 2 R f(0)= 2R |K\cap \mathbb R^{m}| . $$
This finishes the proof. $\hfill\Box $ 

\smallskip

\begin{proposition}
Let $K$ be a symmetric convex body in $\mathbb R^n$, $n\gr 3$ and $2\ls k \ls n-1$. Assume that there exists $\delta\gr 1$ such that for every $F_{1}, F_{2}\in G_{n,k}$,
$$ \frac{ |K\cap F_{1}|}{ |K\cap F_{2}|}  \ls  \delta^{k} . $$
Then 
$$ d_{G}(K, B_{2}^{n}) \ls k\delta^{k} . \leqno(4.10) $$
\end{proposition}

\noindent {\it Proof.}
\noindent Let $\theta_{1}\in S^{n-1}$ be such that $R(K):= h_{K}(\theta_{1})= \rho_{K}(\theta_{1})$. Let $\theta_{2}\in S^{n-1}$ be such that 
$r(K):= h_{K}(\theta_{2})= \rho_{K}(\theta_{2})$. It is enough to show that $\frac{R}{r}\ls k\delta^{k}$. 

\noindent We may assume that $\theta_{1}\neq \theta_{2}$, or else we have nothing to prove. Let $F_{0}:= {\rm span}\{ \theta_{1}, \theta_{2}\} $.
Let $F\in G_{n,k-1}$, $F\perp F_{0}$, $F_{1}:=  {\rm span}\{ \theta_{1},F\}$ and $F_{2}:=  {\rm span}\{ \theta_{2},F\}$. Let $S:= |K\cap F|$. 
Then by Lemma 4.5 we have that 
$$ \frac{2RS}{k} \ls  |K\cap F_{1}| \ls 2R S \ {\rm and} \  \frac{2rS}{k} \ls  |K\cap F_{2}| \ls 2r S \ {\rm or} $$
$$ \frac{R}{kr} \ls \frac{ |K\cap F_{1}|}{ |K\cap F_{2}|} \ls \frac{kR}{r} \leqno(4.11)$$ 
Then $\frac{R}{kr}\ls \delta^{k}$. $\hfill\Box $ 

\smallskip

\medskip

\noindent {\it Proof of Theorem 1.2:} The case $k=1$ is covered by the result of Hensley as noted in the Introduction. So we assume that $n\gr 3$ and $2\ls k \ls n-1$.
Using (4.2) we may assume that $K$ is isotropic. Then we want to show that $d_{G}(I_{k}(K), B_{2}^{n}) \ls k(ck)^{k}$, where $c>0$ absolute constant. 
This follows from Corollary 4.4 and Proposition 4.6. $\hfill\Box $ 

\smallskip

\noindent The dependence on $k$ in Theorem 1.2 is very bad, $c(k)\ls k(ck)^{k}$. This becomes meaningless for $k\gr \frac{\log{n}}{\log\log{n}}$. 
We can give a (slightly) better bound (but still exponential) using certain tools that were developed in order to attack the Hyperplane conjecture. 
For a proof of the following lemma see \cite{MP}.

\begin{lemma}
\noindent Let $K$ be an isotropic convex body in $\mathbb R^n$ and let $1\ls k \ls n-1$. Then for any $F\in G_{n,k}$ there exists a symmetric 
convex body $B$ in $F$ such that, $B:= B(F)$ is also isotropic and  
$$ |K\cap F^{\perp}|^{\frac{1}{k}}  \simeq \frac{L_{B}}{L_{K}} . \leqno(4.12) $$
\end{lemma}

\noindent {\it Second proof of Theorem 1.2:} Again we may assume that $K$ is isotropic.
We will use the best known bound for the isotropic constant due to B. Klartag \cite{Kl} (see also \cite{EMKl}): for every $B$ in $\mathbb R^k$,
$$ L_{B} \ls c k^{\frac{1}{4}} , \leqno(4.13) $$
where $c>0$ is a universal constant. Moreover, it is known (e.g.(\cite{MP}) that $L_{K}\gr L_{B_{2}^{n}}\simeq 1$. So using (1.12) we get that 
for every $F_{1}, F_{2} \in G_{n,k}$,
$$ \frac{ | I_{k}(K) \cap F_{1}|}{  | I_{k}(K) \cap F_{2}|}= \frac{ |K\cap F_{1}^{\perp}|}{ |K\cap F_{2}^{\perp}|} \ls \left(c \frac{L_{B(F_{1})}}
{ L_{B(F_{2})}}\right)^{k} \ls \left(c^{\prime}k\right)^{\frac{k}{4}} . \leqno(4.14) $$
So, using Proposition 4.6 again we get that $d_{G}(K, B_{2}^{n}) \ls k (ck)^{\frac{k}{4}}$. This finishes the proof. $\hfill\Box $ 

\smallskip


\section{Volumetric estimates for $k$-intersection bodies}

\noindent Let $p\neq 0$ and $C$ be a symmetric star body in $\mathbb R^n$. We define
$$ M_{p}(C):= \left(\int_{S^{n-1}} \|\theta\|_{C}^{p} d \sigma(\theta)\right)^{\frac{1}{p}}. \leqno(5.1) $$
Moreover, if $C$ is convex we write
$$ W_{p}(C):= \left( \int_{S^{n-1}} h_{C}^{p}(\theta) d\sigma(\theta) \right)^{\frac{1}{p}} . \leqno(5.2) $$

\noindent Let $K$ be a compact set in $\mathbb R^n$ with $|K|=1$. Let $p>-n, \  p\neq 0$. We define
$$ E_{p}(K) := \left( \int_{K} \|x\|_{2}^{p} d x\right)^{\frac{1}{p}} . \leqno(5.3) $$ 

\noindent We have the following identity (see \cite{Pa4} for a generalization to the case of measures):

\begin{lemma}
\noindent Let $K$ be a symmetric star body in $\mathbb R^n$ of volume $1$ and $1\ls k \ls n-1$. Then 
$$ E_{-k}(K) \left( \int_{G_{n,k}} |K\cap F^{\perp}| d \mu(F) \right)^{\frac{1}{k}} =  E_{-k}(D_{n}) \left( \int_{G_{n,k}} |D_{n}\cap 
F^{\perp}| d \mu(F) \right)^{\frac{1}{k}} . \leqno(5.4) $$
\end{lemma}

\noindent {\it Proof.}
\noindent Indeed, writing in polar coordinates, we get that 
$$ E_{-k}^{-k}(K) = \frac{n\omega_{n}}{ n-k} \int_{G_{n,n-k}} \int_{S_{E}} \frac{d\sigma(\theta)}{\|\theta\|_{K}^{n-k}} d \mu(E) = $$
$$ \frac{n\omega_{n}}{ (n-k)\omega_{n-k}} \int_{G_{n,n-k}} \omega_{n-k}  \int_{S_{E}} \frac{d\sigma(\theta)}{\|\theta\|_{K}^{n-k}} d \mu(E) =  
\frac{n\omega_{n}}{ (n-k)\omega_{n-k}} \int_{G_{n,n-k}}  |K\cap E| d\mu(E) = $$
$$ \frac{n\omega_{n}}{ (n-k)\omega_{n-k}} \int_{G_{n,k}}  |K\cap F^{\perp}| d\mu(F) . $$
Since the same holds also for $D_{n}$, equation (5.4) follows. $\hfill\Box $ 

\smallskip

\noindent We have the following application of the previous lemma and a definition of $k$-intersection bodies. 

\begin{lemma}
\noindent Let $K$ be a symmetric star body in $\mathbb R^n$ of volume $1$ and $1\ls k \ls n-1$. Assume that $I_{k}(K)$ exists. Then 
$$ M_{-k}(I_{k}(K))  \left( \int_{G_{n,k}} |K\cap F^{\perp}| d \mu(F) \right)^{\frac{1}{k}} = \omega_{k}^{\frac{1}{k}} .\leqno(5.5)$$
Moreover,
$$ \frac{M_{-k}(I_{k}(K)) }{ M_{-k}(I_{k}(D_{n}))} = \frac{ E_{-k}(K)}{ E_{-k}(D_{n}) } .\leqno(5.6)  $$ 
\end{lemma}

\noindent {\it Proof.}
\noindent Integrating (4.1) over $G_{n,k}$ we have that 
$$  \int_{G_{n,k}} |I_{k}(K) \cap F  | d\mu(F)  = \int_{G_{n,k}} |K\cap F^{\perp}| d\mu(F) . \leqno(5.7) $$
Moreover, 
$$  \int_{G_{n,k}} |I_{k}(K) \cap F  | d\mu(F) = \int_{G_{n,k}} \omega_{k} \int_{S_{F}} \frac{d\sigma_{F}(\theta)}{ \|\theta\|_{I_{k}(K)}^{k}} d\mu(F)  $$
$$ =\omega_{k} \int_{S^{n-1}}\frac{d\sigma(\theta)}{ \|\theta\|_{I_{k}(K)}^{k}} =  \omega_{k} M_{-k}^{-k} (I_{k}(K)) . $$
So, by (5.7) we get (5.5). (5.6) follows from (5.5) and (5.4). $\hfill\Box $ 

\smallskip

\noindent The following two propositions deal with the behavior of the ratio $ \frac{ E_{p}(K)}{ E_{p}(D_{n}) } $.

\noindent It is not difficult to obtain a lower bound for the quantity $\frac{E_{p}(K)}{E_{p}(D_{n})}$. The main tool is an argument of 
Milman and Pajor (see \cite{MP}).

\begin{proposition}
\noindent Let $K$ a compact set in $\mathbb R^n$ of volume $1$ and let $p>-n, \  p\neq 0$. Then 
$$ E_{p}(K) \gr E_{p} (D_{n}), \leqno(5.8)$$
\noindent with equality if and only if $|K \cap D_{n}|=1$.

\noindent In particular, if $K$ is a star body, then we have equality if and only if $K=D_{n}$. 
\end{proposition}

\noindent{\it Proof. }
\noindent Note that $|K\setminus D_{n}| = | D_{n} \setminus K|$. Also if $x \in K\setminus D_{n}$ then $\|x\|_{2} \gr |B_{2}^{n}|^{-\frac{1}{n}}$ 
while if $x \in  D_{n} \setminus K$, $ \|x\|_{2} \ls |B_{2}^{n}|^{-\frac{1}{n}}$. So, if $p>0$,
$$ E_{p}^{p}(K) = \int_{K} \|x\|_{2}^{p} d x = \int_{K\setminus D_{n}} \|x\|_{2}^{p} d x + \int_{K\cap  D_{n}} \|x\|_{2}^{p} d x \gr $$
$$ \int_{ D_{n} \setminus K} \|x\|_{2}^{p} d x + \int_{K\cap  D_{n}} \|x\|_{2}^{p} d x = E_{2}^{p}( D_{n})  .$$

\noindent It is clear that we have equality if and only if $|K \cap D_{n}|=1.$ We work similarly if $p<0$.
$\hfill\Box $ 

\begin{proposition}
\noindent Let $K$ a star body in $\mathbb R^n$ of volume $1$, $-n < p < q \ls \infty$, $p,q \neq 0$. Then
$$ \frac{E_{p}(K)}{E_{q}(K)} \ls \frac{E_{p}( D_{n})}{E_{q}( D_{n})} , \leqno(5.9) $$
\noindent with equality if and only if $K= D_{n}$. 
\end{proposition}

\noindent {\it Proof.} We follow an argument from \cite{BK}. A simple computation shows that
$$\frac{E_q(D_{n})}{E_p(D_{n})}=
\frac{(\frac{n}{n+q})^{1/q}}{(\frac{n}{n+p})^{1/p}} . \leqno
(5.10) $$

\noindent For every $q>-n$, $q\neq 0$, by integration in polar coordinates we have 

$$ E_{q}^{q}(K)=\omega_n \int_{0}^{\infty} r^{n+q-1}
\sigma\left (\frac{1}{r}K\right )\, d r := \ \int_{0}^{\infty} r^{n+q-1}
g(r) , d r .\leqno (5.11)$$ 

\noindent The function $g(r):=\omega_n\sigma\left( \frac{1}{r} K \right)$ is
non-increasing on $(0,\infty )$ and it is supported in some $[r(K), R(K)]$. If we assume that $K$ is not the Euclidean ball, then $[r(K), R(K)]$ is 
an interval and $g(r)$ can be assumed absolutely
continuous. In this case we can write
$$ g(r)=n \int_{r}^{\infty} \frac{\rho (s)}{s^n} d s ,\,\quad (r>0)\leqno
(5.12)$$ 
for some non-negative function $\rho $ on $(0, \infty)$.
Then, again by integration in polar coordinates,
$$ 1=|K|= \int_{0}^{\infty} r^{n-1}g(r)\,d r = n\int_0^{\infty }\int_{0<r<s} r^{n-1}
\frac{\rho (s)}{s^n}d r\,d s = \int_{0}^{\infty}\rho (s)\,d s.$$
Hence, $\rho $ represents a probability density of a
positive random variable, say, $\xi$. Then (5.11) becomes

$$E_{q}^{q}(K)=\int_0^{\infty }r^{q+n-1}g(r)\,d r = \frac{n}{n+q}
\int_{0}^{\infty} s^q \rho (s)\,d s =\frac{n}{n+q} {\mathbb E}
(\xi^q) .$$
Applying H\"{o}lder's inequality for
$-n < p <  q \ls \infty$, we see that
$$\left({\mathbb E} (\xi^q)\right)^{1/q} > \left({\mathbb E}(\xi^p)
\right)^{1/p} \ . \leqno(5.13)$$
Note that, since $\xi$ is a non-zero random variable with an absolute continuous density, there is no equality case in (5.13) (see \cite[Th.188]{HLP}).
So,
$$\frac{E_{q}(K)}{E_{p}(K)}=\frac{\left(\frac{n}{n+q} {\mathbb E} (\xi^q)\right)^{1/q}}
{\left(\frac{n}{n+p} {\mathbb E} (\xi^p)\right)^{1/p}} >
\frac{\left(\frac{n}{n+q}\right)^{1/q}}{\left(\frac{n}{n+p}\right)^{1/p}}=
\frac{E_q(D_{n})}{E_p(D_{n})}, $$
as claimed. 
$\hfill\Box $

\smallskip

\noindent We will use the following immediate application of H\"older's inequality. 

\begin{lemma}
\noindent Let $C$ be a star symmetric body in $\mathbb R^n$ and $p\ls q$, $p,q\neq 0$. Then 
$$ M_{p}(C) \ls M_{q}(C) , \leqno(5.14) $$
with equality if and only if $C=aB_{2}^{n}$ for some $a>0$. 
\end{lemma}

\noindent Moreover writing the volume of $C$ in polar coordinates  we get that
$$ M_{-n}(C) = \left(\frac{ |B_{2}^{n}|} {|C|}\right)^{\frac{1}{n}}. \leqno(5.15) $$

\noindent Results of Lewis \cite{Lew}, Figiel and Tomczak-Jaegermann \cite{FT}, Pisier \cite{Pisier1} establish the following ``reverse Uryson" inequality:

\begin{proposition}
\noindent Let $C$ be a symmetric convex body in $\mathbb R^n$. Then there exists $T\in SL_{n}$ such that 
$$ W(C_{1}) \ls c \sqrt{n} |C_{1}|^{\frac{1}{n}} \log{\left(1+d_{BM}(C, B_{2}^{n})\right)} , \leqno(5.16)$$
where $c>0$ is an absolute constant and $C_{1}= TC$.
\end{proposition}

\medskip

\noindent {\it Proof of Theorem 1.3:} Using Lemma 5.5 and (5.15) we have that 
$$  \frac{M_{-k}(I_{k}(K)) }{ M_{-k}(I_{k}(D_{n}))}  \gr \left(\frac{ |I_{k}(D_{n})|} {|I_{k}(K)|} \right)^{\frac{1}{n}} . \leqno(5.17)$$
Using (4.2), we may assume that $K$ is isotropic. Then $E_{2}(K)= \sqrt{n}L_{K}$. Then, by (5.6), and Proposition 5.4,
$$ \frac{M_{-k}(I_{k}(K)) }{ M_{-k}(I_{k}(D_{n}))} = \frac{ E_{-k}(K)}{ E_{-k}(D_{n}) } \ls  \frac{ E_{2}(K)}{ E_{2}(D_{n}) }  = \frac{L_{K}}{L_{D_{n}}} . \leqno(5.18) $$
So, (5.17) and (5.18) imply equation (1.3). 

\smallskip

\noindent We now assume that $I_{k}(K)$ is convex. We consider $I_{k}^{\circ}(K)$ to be in the position described in Proposition 5.6. (Again by using (4.2).)
Then by Lemmas 5.5 and 5.6, 
$$ M_{-k}(I_{k}(K)) \ls M_{1}(I_{k}(K)) = W (I_{k}^{\circ}(K)) \ls $$
$$ c \sqrt{n} | |I_{k}^{\circ}(K) |^{\frac{1}{n}} \log{\left(1+d_{BM}(I_{k}(K), B_{2}^{n})\right)} \ls \frac{ c^{\prime}\log{\left(1+d_{BM}(I_{k}(K), B_{2}^{n})\right)}  }
{ \sqrt{n} |I_{k}(K)|^{\frac{1}{n}} } , $$
using also Santal\'o inequality.  This implies that 
$$  \frac{M_{-k}(I_{k}(K)) }{ M_{-k}(I_{k}(D_{n}))} \ls c  \left(\frac{ |I_{k}(D_{n})|} {|I_{k}(K)|} \right)^{\frac{1}{n}}\log{\left(1+d_{BM}(I_{k}(K), B_{2}^{n})\right)}   . \leqno(5.19)$$
So again by (5.6) and Proposition 5.3, we have that 
$$ \frac{M_{-k}(I_{k}(K)) }{ M_{-k}(I_{k}(D_{n}))} = \frac{ E_{-k}(K)}{ E_{-k}(D_{n}) } \gr 1 . \leqno(5.20) $$
So, equations (5.19) and (5.20) imply (1.4). To conclude we use the fact that $d_{BM}(K, B_{2}^{n}) \ls \sqrt{n}$ and Theorem 1.2.  $\hfill\Box $

\bigskip


\section{Concluding remarks}

\noindent Theorem 1.1 indicates that the classes ${\mathcal{BP}}_{k}^{n}$ increase in a ``canonical" way. It is not clear to us if the reverse inequality (up to the logarithmic term) holds true even for the Banach-Mazur distance. We pose this as a 

\smallskip

\noindent {\bf Question 1:} Is it true that for every $n$ and every $k,$ $1\ls k \ls n-1$, there exists a convex symmetric body $K$ in $\mathbb R^n$  such that
$$ d_{BM} (K, {\mathcal{BP}}_{k}^{n}) \gr c\sqrt{\frac{n}{k}} , $$
where $c>0$ is a universal constant?  

\smallskip

\noindent We can show that this is true in the case $k=1$:

\begin{proposition}
\noindent There exists $c>0$ such that for every $n\gr 1$, 
$$ d_{BM}( B_{\infty}^{n}, {\mathcal{I}}_{1}^{n}) \gr c \sqrt{n} . $$
\end{proposition}

\noindent The proof of the latter Proposition depends on the following fact: all convex intersection bodies have bounded volume ratio. Let ${\mathcal{VR}}_{n}(a)$ be the class of symmetric convex bodies in $\mathbb R^n$ with ${\rm v.r. }(K) := \inf\{ \left( \frac{ |TK|}{B_{2}^{n}}\right)^{\frac{1}{n}} : T\in GL_{n}\} \ls a $. In this notation we have the following

\begin{proposition}
\noindent There exists $c>0$ such that for every $n\gr 1$ every convex intersection body in $\R^n$
belongs to the class ${\mathcal{VR}}_{n}(c) .$
\end{proposition}

\noindent {\it Proof.}
\noindent The proof is simply a combination of certain known results. Let $K\in {\mathcal{I}}_{1}^{n}$ and convex, and let $X$ the $n$-dimensional Banach space with norm $\|\cdot\|:= \| \cdot \|_{K}$. Then it has been shown in {\cite{KK2}} that $ X$ embeds isomorphically to $L_{1/2}$. It is known (e.g. \cite{Kalton1}) that every Banach subspace of $L_{1/2}$ has (Rademacher) cotype 2. Next, by a result of Bourgain and Milman (\cite{BM2}), every finite dimensional subspace with bounded cotype $2$ constant, has the property that its unit ball has bounded volume ratio. This finishes the proof.  $\hfill\Box $ 

\smallskip

\noindent {\it Proof of Proposition 6.1:} We will show that $ d_{BM} ( B_{\infty}^{n} , {\mathcal{VR}}_{n}(a)) \gr c(a) \sqrt{n}$. Then by Proposition 6.2 the proof would be complete. For simplicity we assume that $n=4 k$ for some $k\in \mathbb N$. The general case follows easily.  Let $C\in{\mathcal{VR}}_{n}(a)$ such that $  d_{BM} ( B_{\infty}^{n} , {\mathcal{VR}}_{n}(a)) = d_{G}(B_{\infty}^{n}, C)=: d$. Then we have that $ C\subseteq B_{\infty}^{n} \subseteq d C$. By a well known generalization of Kashin's theorem (\cite{ST}) there exists $ F\in G_{n, \frac{n}{2}}$ and an ellipsoid ${\mathcal{E}}$ such that $c_{1}(a){\mathcal{E}} \subseteq C\cap F \subseteq c_{2}(a){\mathcal{E}}$. Moreover it is known (see e.g. \cite{Zg}) that for any ellipsoid in $\mathbb R^m$ there exists $ E\in G_{m, \frac{m}{2}}$ and some $r>0$ such that $ {\mathcal{E}}\cap E= r B_{E}$, where $B_{E}$ is the Euclidean ball of $E$. So we get that there exists $E\in G_{n, \frac{n}{4}}$ and $r_{1}:=r_{1}(a)>0$, $r_{2}:=r_{2}(a)$, such that 
$$ r_{1} B_{E} \subseteq B_{\infty}^{n}\cap E \subseteq d r_{2} B_{E} . $$
Hence, it is enough to show that $d_{G}( B_{\infty}^{n}\cap E, B_{E}) \gr c \sqrt{n}$, for every $E\in G_{n, \frac{n}{4}}$. Considering the polar body of $B_{\infty}^{n}\cap E$ it is enough to show that the convex hull of at most $8n$ points in $\mathbb R^n$ has geometric distance from the Euclidean ball at least $c\sqrt{n}$. Let $N\gr n+1$, $v_{1}, \cdots , v_{N}\in B_{2}^{n}$ and $\|v_{1}\|_{2}=1$. Let $K:= {\rm co}\{ v_{1}, \cdots, v_{N}\}$. Note that $R(K)=1$. Now it is enough to show that if $N=8n$ then $r_{K}:= \min_{\theta_{S^{n-1}}} \rho_{K}(\theta) \ls \frac{c}{\sqrt{n}}$.  But (see e.g. \cite{BF1}, \cite{BF2}, \cite{CP}, \cite{Glu}) one has that $|K|^{\frac{1}{n}} \ls c^{\prime} \frac{ \sqrt{\log{\frac{eN}{n}}}}{ n} $. Writing the volume of $K$ in polar coordinates we get that there exists at least one $\theta \in S^{n-1}$ such that $ \rho_{K}(\theta) \ls c^{\prime\prime} \sqrt{\frac{\log{\frac{eN}{n}}}{n}}$. We complete the proof by choosing $N=8n$ . $\hfill\Box $ 

\medskip

\noindent The estimate in Theorem 1.2 is exponential with respect to $k$. Even if we assume that the Hyperplane conjecture has a positive answer the existing proof would still give an estimate exponential with respect to $k$. We believe that a better estimate (polynomial) must be true. Having in mind equation (4.6), we pose the following question: 

\smallskip

\noindent {\bf Question 2:} Is it true that if $K$ is symmetric and convex in $\mathbb R^n$, $1\ls k \ls n$ and $I_{k}(K)$ exists and it is convex then $I_{k}(K)$ is isomorphic to $Z_{k}^{\circ}(K)$? 

\smallskip

\noindent Note that a positive answer to the previous question would easily imply a linear in $k$ estimate in Theorem 1.2. 

\medskip

\noindent The second conclusion of Theorem 1.3. can be viewed as a generalization of the classical Busemann inequality (see Introduction). 
We don't know if the assumption that $I_{k}(K)$ is convex is necessary in Theorem 1.3 and whether the estimate can be replaced by $1$: 

\smallskip

\noindent {\bf Question 3:} Is it true that if $K$ is a symmetric star body of volume $1$, $1\ls k \ls n-1$ and if $I_{k}(K)$ exists, then 
$$ |I_{k}(K)| \ls | I_{k}(D_{n})| \ ? $$  

\bigskip

{\bf Acknowledgments:} The first name author wish to thank the US National Science Foundation for support through grants DMS-0652571 and DMS-1001234, and the Max Planck Institute for Mathematics for support and hospitality during his stay in Spring 2011. The second name author was partially supported by NSF 0906051. Part of this work was carried out when the third name author was visiting the Mathematics Department of Texas A$\&$ M which she thanks for its hospitality.

\footnotesize
\bibliographystyle{amsplain}

\end{document}